\documentclass[11pt,a4paper]{amsart}
\usepackage{graphicx}
\usepackage{amsmath,amssymb}
\usepackage{epstopdf}
\usepackage{pgf, pgfnodes}
\usepackage{tikz}
\usetikzlibrary{calc}
\usepackage{url}
\usepackage{geometry}

\newtheorem{prop}{Proposition}

\def\sqr#1#2{{\vbox{\hrule height.#2pt
    \hbox{\vrule width.#2pt height#1pt \kern#1pt
        \vrule width.#2pt}\hrule height.#2pt}}}
\def\eqed{\sqr53}
\def\qed{%
    \ifmmode\eqno\eqed
    \else\nobreak\ \hfill\eqed\medbreak\fi}

\newcommand{\ser}{\mathbin{\bowtie}}
\newcommand{\pll}{\mathbin{\|}}

\title{The Beraha number $B_{10}$ is a chromatic root}
\author{Gordon Royle}

\begin{document}
\maketitle
\begin{abstract}
This note exhibits graphs whose chromatic polynomials have the Beraha number $B_{10} = (5 + \sqrt{5})/2$ as
a root. It was previously known that no other non-integer Beraha number is a chromatic root, and so 
these examples complete the determination of precisely which Beraha numbers can be chromatic roots. 
\end{abstract}

\thispagestyle{empty}

\section{Introduction}

Given a graph $G$, its \emph{chromatic polynomial} $P_G$ is the function defined by the property that $P_G(q)$ is the number of proper $q$-colourings of $G$ whenever $q$ is a non-negative integer. As its name suggests, $P_G$ turns out be a polynomial in $q$ and so it is possible to evaluate the function at arbitrary real and complex values of $q$, regardless of whether these evaluations have any combinatorial significance. In particular, it is possible to find the roots of this polynomial, both real and complex, and there is an extensive body of work relating graph-theoretical properties to the location of these \emph{chromatic roots}. The fundamental connection between complex chromatic roots and phase transitions in the $q$-state Potts model means that much of these results appear in the statistical physics literature.

The fundamental results regarding which real numbers can be chromatic roots (of any graph) can be summarised in a 
single sentence: There are no chromatic roots in the real intervals $(-\infty,0)$, $(0,1)$ or $(1,32/27]$ (Jackson \cite{MR1264037}), while chromatic roots are dense in the interval $(32/27,\infty)$ (Thomassen \cite{MR1483433}). (There are more fine-grained results for specific classes of graphs, such as planar graphs, but we do not need these now.)

The \emph{Beraha numbers} are an infinite sequence of numbers $\{B_n\}_{n=1}^{\infty}$ given by
\begin{equation}
B_n = 2 + 2 \cos \left( \frac{2 \pi}{n} \right).
\end{equation}
The sequence starts 
\[
B_1 = 4,\ 0,\ 1,\ 2, \varphi+1,\ 3,\  3.2469796,\ 2 + \sqrt{2},\ 3.5320889,\ \varphi+2,\ \ldots
\]
where the values for $B_7$ and $B_9$ are rounded to 8 significant figures, and $\varphi = (1+\sqrt{5})/2$ is the Golden Ratio. The sequence was first identified by Beraha, based on empirical observations that certain real numbers seem to play a particularly important role in the study of chromatic roots of planar graphs, especially as limit points of complex zeros of the various lattices studied by statistical physicists. His conjecture that there are chromatic roots of planar triangulations arbitrarily close to each Beraha number is still not completely resolved.

One of the more remarkable facts involving some of the Beraha numbers is Tutte's famous ``Golden Identity'' (Tutte \cite{MR0272676}) relating the values of the chromatic polynomial of a planar triangulation at $B_{5}$ and $B_{10}$. The Golden Identity asserts that if $T$ is a planar triangulation on $n$ vertices, then
\begin{equation}
P_T(\varphi+2) = (\varphi+2) \varphi^{3n-10} P_T(\varphi+1)^2,
\end{equation}
with the implication that for any planar triangulation $T$, either both or neither $B_5$ and $B_{10}$ are chromatic roots of $T$. In 1970, Tutte \cite{MR0272676}  already knew that  $P_T(B_5) \not= 0$ for any planar triangulation $T$ and therefore that $P_T(B_{10})$ is not only non-zero, but strictly positive. Recently Perrett and Thomassen (personal communication) have shown that the same conclusion holds for all planar graphs, not just triangulations.  

Salas and Sokal \cite{MR1853428} observed that not only is $B_5$ not the chromatic root of any planar triangulation, but in fact it is not the chromatic
root of any graph whatsoever. This follows because the minimal polynomial of $B_5$ is $x^2 - 3x +1$ and so any integer polynomial with $B_5 = (3 + \sqrt{5})/2$ as a root, must also have its algebraic conjugate (the other root of the minimal polynomial) $B_5^* = (3 - \sqrt{5})/2$ as a root. As $B_5^* \approx 0.38196601$ lies in the forbidden interval $(0,1)$ it follows that neither $B_5$ nor
$B_5^*$ are chromatic roots of any graph.  In a similar fashion, they showed that no non-integer Beraha number is the chromatic root of any graph, with the possible exception of $B_{10} = \varphi+2$. The minimal polynomial of $B_{10}$ is $x^2 - 5 x + 5$ and so the algebraic conjugate $B_{10}^* \approx 1.381966011$ which is not \emph{a priori} forbidden.  They concluded that \emph{``The exceptional case $n=10$ is very curious''}.

In the rest of this short note, we resolve this exceptional case by exhibiting a graph with $B_{10}$ as a chromatic root, thereby completing the determination of exactly which Beraha numbers can be chromatic roots.

\begin{prop}
The only Beraha numbers that are chromatic roots are the integer Beraha numbers $\{B_1, B_2, B_3, B_4, B_6\}$ and $B_{10}$.
\end{prop}

\section{The Graphs}

\begin{figure}
\begin{center}
\begin{tikzpicture}[rotate=90,scale=0.95]
\tikzstyle{vertex}=[circle, draw=black, fill=black!25!white, inner sep = 0.55mm]
\tikzstyle{svertex}=[circle, draw=black, fill=black!25!white, inner sep = 0.45mm]
\node [vertex] (v0) at (30:2cm) {};
\node [vertex] (v1) at (90:2cm) {};
\node [vertex] (v2) at (135:2cm) {};
\node [vertex] (v3) at (225:2cm) {};
\node [vertex] (v4) at (270:2cm) {};
\node [vertex] (v5) at (330:2cm) {};

\node [vertex] (v6) at (180:2.1cm) {};

\node [svertex] (v7) at (157:1.8cm) {};
\node [svertex] (v8) at (157:2.2cm) {};

\node [svertex] (v9) at (203:1.8cm) {};
\node [svertex] (v10) at (203:2.2cm) {};

\draw (v0)--(v1)--(v2)--(v0);
\draw (v0)--(v3);
\draw (v0)--(v4);
\draw (v0)--(v5);
\draw (v3)--(v4)--(v5)--(v3);
\draw (v1)--(v3);
\draw (v1)--(v4);
\draw (v1)--(v5);
\draw (v2)--(v4);
\draw (v2)--(v5);

\draw (v2)--(v7)--(v6)--(v8)--(v2);
\draw (v3)--(v9)--(v6)--(v10)--(v3);

\end{tikzpicture}
\hspace{2cm}
\begin{tikzpicture}
\tikzstyle{vertex}=[circle, draw=black, fill=black!25!white, inner sep = 0.55mm]
\tikzstyle{svertex}=[circle, draw=black, fill=black!25!white, inner sep = 0.45mm]
\node [vertex] (v2) at (90:2cm) {};
\node [vertex] (v3) at (162:2cm) {};
\node [vertex] (v4) at (234:2cm) {};
\node [vertex] (v0) at (306:2cm) {};
\node [vertex] (v1) at (18:2cm) {};

\node [svertex] (v5) at ($(v4) + (0.75,0.35)$) {};
\node [svertex] (v6) at ($(v4) + (0.75,0)$) {};
\node [svertex] (v7) at ($(v4) + (0.75,-0.35)$) {};

\node [svertex] (v8) at ($(v0) - (0.75,0.35)$) {};
\node [svertex] (v9) at ($(v0) - (0.75,0)$) {};
\node [svertex] (v10) at ($(v0) - (0.75,-0.35)$) {};

\draw (v0)--(v1);
\draw (v0)--(v2);
\draw [bend left = 15] (v4) to (v1);
\draw [bend left = 15] (v3) to (v0);
\draw (v4)--(v2);
\draw (v4)--(v3);

\draw (v1)--(v2)--(v3)--(v1);

\draw (v4)--(v5);
\draw (v4)--(v6);
\draw (v4)--(v7);

\draw (v0)--(v8);
\draw (v0)--(v9);
\draw (v0)--(v10);

\draw (v10)--(v5);
\draw (v10)--(v6);

\draw (v7)--(v8);
\draw (v7)--(v9);

\end{tikzpicture}
\end{center} 
\caption{Two graphs $G_1$ (left) and $G_2$ (right) with $B_{10}$ as a chromatic root}
\label{smallest}
\end{figure}
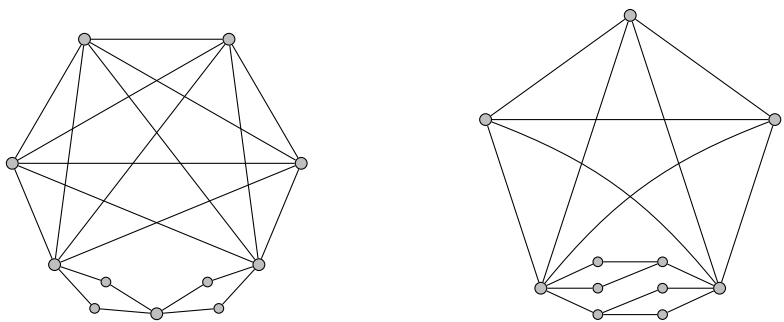

Two graphs with $B_{10}$ as a chromatic root are shown in Figure~\ref{smallest}. Each example is obtained by replacing a single edge of a suitable complete graph (either $K_6$ or $K_5$) with a small ``gadget'' as shown in Figure~\ref{smallest}. Each of the  gadgets is a small $2$-terminal series-parallel graph (see  Royle and Sokal \cite{MR3416855} for an extensive discussion of $2$-terminal series-parallel graphs). The factored chromatic polynomials of these graphs are 
{\small 
\begin{align*}
P_{G_1}(q) &= q (q-1) (q-2) (q-3)^2 (q-4) \underbrace{\left(q^2-5 q+5\right)}_{\text{Min. poly of $B_{10}$}} \left(q^3-4 q^2+8 q-7\right)\\
P_{G_2}(q) &= q (q-1) (q-2) (q-3)  \underbrace{\left(q^2-5 q+5\right)}_{\text{Min. poly of $B_{10}$}}\left(q^5-8 q^4+30 q^3-63 q^2+73 q-36\right),
\end{align*}}
with the factor $q^2-5q+5$ highlighted.

Of course it is easy to verify these chromatic polynomials directly by computer, but it is nonetheless interesting to 
investigate the structure of these examples and derive the result directly and in a form that can easily be
verified without using chromatic polynomial software. This is possible due to the particular convenience of chromatic polynomial
calculations (in fact, even Tutte polynomial calculations) in $2$-terminal graphs (and particularly $2$-terminal series-parallel graphs). All that is required is to maintain some auxiliary information regarding how the terminals are 
coloured and this information can then be propagated through the operations of series and parallel connections. 

There are a variety of ways of doing this, but for our current purposes we will express the chromatic polynomial of a $2$-terminal graph as the sum of the $q$-colourings that colour the terminals with the {\em same} colour and those that colour the 
terminals with {\em different} colours. So for a 2-terminal graph $G$ with terminals $s$, $t$, we let $S_G(q)$ denote the number of $q$-colourings of $G$ where $s$ and $t$ are coloured the {\bf S}ame, and $D_G(q)$ the number of $q$-colourings of $G$ where they are coloured {\bf D}ifferent.  Then it is clear that $P_G(q) = S_G(q) + D_G(q)$ and that
\begin{align*}
S_G(q) &= P_{G/st}(q) \\
D_G(q) &= P_{G+st}(q) 
\end{align*}
where $G/st$ denotes the graph obtained by \emph{merging} $s$ and $t$ (creating a loop if $s$ and $t$ are adjacent), and $G+st$ the graph obtained by \emph{joining} $s$ and $t$. 

So now suppose that $G$ and $H$ are 2-terminal graphs and that we know $S_G$, $S_H$, $D_G$ and $D_H$. Then
 we can determine the  same values for the parallel connection $G \pll H$ and the series connection $G \ser H$ as 
 follows: 
\begin{align*}
S_{G \| H} (q) &= \frac{S_G(q) S_H(q)}{q}, \\
D_{G \| H} (q) &= \frac{D_G(q) D_H(q)}{q(q-1)},\\
S_{G \ser H} (q) &= \frac{S_G (q) S_H(q)}{q} + \frac{D_G(q) D_H(q) }{ q(q-1) }, \\
D_{G \ser H} (q) &= \frac{(q-2) D_G(q) D_H(q)}{q(q-1)} + \frac{D_G(q) S_H(q)}{q} + \frac{S_G(q) D_H(q)}{q} .
\end{align*}
These expressions are calculated by taking the Cartesian product of the $q$-colourings of $G$ and $H$, dividing by some function of $q$ to ensure that we count only those pairs of colourings that match on the common vertices, and then assigning the resulting colourings to either $S$ or $D$ depending on whether the new terminals are coloured the same or different.

Now we will use these expressions to calculate the chromatic polynomial for $G_1$, which in the notation of 
Royle \& Sokal \cite{MR3416855} can be constructed as follows:
\begin{align*}
W & =  (K_2 \ser K_2) \| (K_2 \ser K_2), \\
G_1 &= (K_6 \backslash e) \pll   (W \ser W),
\end{align*}
where $W$ is the $4$-cycle $C_4$ with nonadjacent vertices as terminals and $K_6 \backslash e$ is viewed as a 2-terminal graph by declaring the former end points of $e$ to be the two terminals.

\newcommand\ff[1]{(q)_#1}

So if we let $p_G(q) = [S_G(q), D_G(q)]$ denote the ``partitioned chromatic polynomial'' of a 2-terminal graph $G$, and let $\ff{n} = q(q-1)\ldots(q-n+1)$ denote the falling factorial (the chromatic polynomial of $K_n$), then we have
\begin{align*}
p_{K_2}(q) &=[ 0, \ff{2} ], \\
p_{K_2 \ser K_2} (q) &= [\ff{2}, \ff{3}], \\
p_{W} (q) &= [ \ff{2} (q-1), \ff{3} (q-2)] ,  \\
p_{W \ser W}  (q) &= [\ff{2} \left(q^4-7 q^3+21 q^2-29 q+15\right), \ff{3} (q-2)  \left(q^3-4 q^2+8 q-6\right)], \\
p_{K_6\backslash e}(q) &= [\ff{5}, \ff{6}]. 
\end{align*}

Finally, combining these last two expressions in the appropriate fashion yields the chromatic polynomial given above. A similar but slightly more complicated derivation yields the chromatic polynomial of $G_2$.

\section{Conclusion}

In some sense this is a slightly disappointing result --- while $B_{10}$ is very significant in the study of chromatic roots of planar graphs, it does not appear to have to any particular significance for general graphs. It seems that the only technique that we have to show that a particular number is not a chromatic root is to find an algebraic conjugate of the number lying in one of the forbidden regions, and until now $B_{10}$ was viewed as a promising candidate for a non-chromatic root with no forbidden conjugates. In fact, Cameron \& Morgan \cite{1610.00424} wonder ``\emph{Perhaps resolving the question whether $\alpha+3$ [i.e. $B_{10}$] is a chromatic root would help with this}", to which we must now respond resoundingly in the negative!

\bibliographystyle{acm2url}
\bibliography{../gordonmaster.bib}

 \end{document}